

\magnification=\magstep1
\voffset=1.00truein
\settabs 18 \columns
\hoffset=1.00truein
\baselineskip=17 pt

\def\s{\smallskip}

\def\b{\bigskip}

\def\bbb{\bigskip\bigskip\bigskip}

\def\sqr#1#2{{\vcenter{\vbox{\hrule height.#2pt
 \hbox{\vrule width.#2pt height#1pt \kern#1pt
 \vrule width.#2pt} \hrule height.#2pt}}}}

\def\perp{\hbox{${\kern+.25em{\bigcirc}
\kern-.85em\bot\kern+.85em\kern-.25em}$}}
\def\lsim{\;raise0.3ex\hbox{$>$\kern-0.75em\raise-1.1ex\hbox{$\sim$}}\;}
\def\gsim{\;raise0.3ex\hbox{$>$\kerm-0.75em\raise-1.1ex\hbox{$\sim$}}\;}
\def\no{\noindent}

\def\ce{\centerline}
\def\ve{\vfill\eject}
\def\rdots{\mathinner{\mkern1mu\raise1pt\vbox{\kern7pt\hbox{.}}\mkern2mu
 \raise4pt\hbox{.}\mkern2mu\raise7pt\hbox{.}\mkern1mu}}

\def\e e{$e^+ e^-$ }


\rightline{UCLA/95/TEP/4}
\rightline{May 1995}
\bbb
\ce {\bf $q$-GAUGE THEORY}
\b
\ce {Robert J. Finkelstein}
\ce {\it Department of Physics}
\ce {\it University of California, Los Angeles, CA 90095-1547}
\bbb
\no {\bf Abstract.}  We examine some issues that arise in the
$q$-deformation of a gauge theory.  If the deformation is
carried out by replacing the equal time commutators of free fields
by the corresponding $q$-commutators, the resulting propagators are not
very much different from those of the undeformed theory as long as one is
dealing with weak fields; but the theory still violates causality.
If one postulates
a $q$-deformed $S$ matrix, the corresponding $q$-causal commutator
has 2 poles of different strength and the result again amounts to
a deformation of the Lorentz group.
\ve
\line {\bf 1. Introduction.  \hfil}
\s

A quantum field theory rests jointly on the underlying geometry
and the dynamical laws operating in this geometry.  The particular
way the structure of the full theory is distributed between these
two foundations is somewhat arbitrary.  In recent years,
beginning with the discoveries of supersymmetry and supergravity,
the emphasis has been on the investigation of different geometries.
Here we shall pursue the alternative option of modifying the
dynamics and more specifically of modifying the quantization
procedure.  This alternative touches on another open and less explored
subject, the rigidity of quantum mechanics.

There is an interesting deformation of quantum mechanics obtained by
altering the Dirac prescription so as to replace classical Poisson
brackets of dynamically conjugate variables by the $q$-commutators
of the corresponding quantum observables.  As far as this proposal
has been investigated there seems to be no obstruction to the
formulation of a $q$-quantum mechanics for finite systems;$^1$ and if
$q$ is close enough to unity to be compatible with present
experiment, these theories may even be regarded as realistic.  One
may try to implement the same idea in field theory by working with fields
and their dynamically conjugate fields.
This is the first point we shall investigate here.
There is an arbitrariness in this formulation
however, and depending on which questions are asked, the resulting theory
may differ either little or greatly from the $q=1$ formalism.
We shall also examine a second approach in which one postulates
a ``$q$-time ordering" of the $S$ matrix.
In the latter case there is a violation of special relativity.  This
is an example in which the altered dynamics is
inconsistent with the original geometry and therefore is not independent
of an explicit $q$-deformation of that geometry---an avenue that has
been much explored.$^2$

The effect of replacing the usual
commutators (or anticommutators) of the field oscillators by
$q$-commutators is to replace the occupation number $n$ by the
corresponding basic number

$$\langle n\rangle = {q^n-1\over q-1}~.$$

\no Then one might anticipate that the Einstein relation leading to
the Planck law, namely:

$${n+1\over n} = e^{h\nu/kT}$$

\no would be replaced by

$${\langle n+1\rangle\over \langle n\rangle}=e^{h\nu_q/kT}$$

\no which is equivalent to

$$\langle n\rangle={1\over e^{h\nu_q/kT}-q}~.$$

This last relation is in fact correct in both the E.B. ($q=1$) and
F.D. ($q=-1$) limits.

It is perhaps of interest to examine the $q$-formulation of other
basic relations holding for the two kinds of statistics.  We shall
here discuss the perturbation sector of $q$-electrodynamics.
\vskip.5cm

\line {\bf 2. Quantization. \hfil}
\s

Quantization may be imposed via the field oscillators as
follows.  Denote the expansion of an arbitrary field by

$$\psi_\alpha(x)=\sum_\rho\bigl[a(\rho)f_\alpha(\rho,x)
  +\bar b(\rho) g_\alpha(\rho,x)\bigr] \eqno(2.1)$$

\no where

$$\sum_\rho = \sum_r \int d\vec p \eqno(2.2)$$
$$f_\alpha(\rho,x)=\biggl({1\over 2\pi}\biggr)^{3/2}
  {u_\alpha(\vec p,r)\over (2p_o)^{1/2}}~e^{-ipx} \eqno(2.3)$$
$$g_\alpha(\rho,x)=\biggl({1\over 2\pi}\biggr)^{3/2}
  {v_\alpha(\vec p,r)\over (2p_o)^{1/2}}~e^{ipx}~. \eqno(2.4)$$

\no Here $a(\bar a)$ and $b(\bar b)$ are the absorption (emission)
operators of particles and antiparticles respectively.  The $\rho$
sum is an integration over momentum and a sum over spin.  The particle
and antiparticle parts of the sum are related by complex conjugation
of the exponentials and by charge conjugation of the spin dependent
functions, according to the following relations:

$$\eqalignno{v(\rho,\vec p)&= C\bar u^T=C(\gamma^o)^T u^* & (2.5)\cr
  u(\rho,\vec p) &= C\bar v^T=C(\gamma^o)^T v^*~. & (2.6)\cr}$$

The quantization of the oscillators may be described by the following
equations:

$$\eqalignno{\bigl(a(\rho),\bar a(\rho^\prime)\bigr)_q &=
  \delta(\rho,\rho^\prime) & (2.7)\cr
  \bigl(b(\rho),\bar b(\rho^\prime)\bigr)_q &=
  \delta(\rho,\rho^\prime) & (2.8)\cr
  \bigl(a(\rho),b(\rho^\prime)\bigr)_q &=0 & (2.9)\cr
  \bigl(\bar b(\rho),\bar a(\rho^\prime)\bigr)_q &= 0 & (2.10)\cr}$$

Introduce the operator ${\cal{C}}$ in Hilbert space which takes a particle
state into an antiparticle state:

$${\cal{C}}~\bar a(\rho)\Psi_o = \epsilon^*\bar b(\rho)
  \Psi_o \eqno(2.11)$$

\no where $|\epsilon|=1$ and $\Psi_o$ is the vacuum state for
which we assume

$${\cal{C}}~\Psi_o=\Psi_o~. \eqno(2.12)$$

\no Then

$$\eqalignno{{\cal{C}}~\bar a~{\cal{C}}^{-1} &= \epsilon^*\bar b & (2.13) \cr
  {\cal{C}}~a~{\cal{C}}^{-1} &= \epsilon b & (2.14)\cr}$$

\no if we also assume

$$\bar {\cal{C}}={\cal{C}}^{-1}~.\eqno(2.15)$$

\no Then

$${\cal{C}}\bigl(a(\rho)\bar a(\rho^\prime)-q\bar a(\rho^\prime)
  a(\rho)\bigr)~{\cal{C}}^{-1}=\delta(\rho,\rho^\prime)  \eqno(2.16)$$

\no implies

$$b(\rho)\bar b(\rho^\prime)-q\bar b(\rho^\prime)
  b(\rho)=\delta(\rho,\rho^\prime)~. \eqno(2.17)$$

\no If we try to add the following $q$-commutators,

$$\eqalignno{\bigl(a(\rho),a(\rho^\prime)\bigr)_q &= 0 & (2.18)\cr
  \bigl(b(\rho),b(\rho^\prime)\bigr)_q &= 0 & (2.19)\cr}$$

\no it is clear that in these relations the only permitted value
of $q$ is $\pm 1$.  We take $q=+1$ and $-1$ for ``bose" and ``fermi"
particles respectively, in (2.18) and (2.19).

\vskip.5cm

\line {\bf 3. Expansion of the $S$-Matrix. \hfil}
\s

In expanding the $S$-matrix one encounters the causal propagator,
$\Delta_F$.  Although the oscillators are now $q$-quantized, it is
still easy to calculate $\Delta_F$ in the absence of background
fields; for the commutator may be expressed in terms of the $q$-commutator
as follows:

$$\bigl(a(\rho),\bar a(\rho^\prime)\bigr)=\bigl(a(\rho),
  \bar a(\rho^\prime)\bigr)_q+(q-1)\bar aa \eqno(3.1)$$

\no and therefore the vacuum expectation value is

$$\langle 0|\bigl(a(\rho),\bar a(\rho^\prime)\bigr)|0\rangle=
  \langle 0|\bigl(a(\rho),\bar a(\rho^\prime)\bigr)_q|0\rangle=
  \delta(\rho,\rho^\prime)~. \eqno(3.2)$$

\no Likewise the anticommutator is

$$\bigl\{a(\rho),\bar a(\rho^\prime)\bigr\}=
  \bigl(a(\rho),\bar a(\rho^\prime)\bigr)_q+(q+1)\bar aa \eqno(3.3)$$

\no and

$$\langle 0|\bigl\{a(\rho),\bar a(\rho^\prime)\bigr\}|0\rangle=
  \langle 0|\bigl(a(\rho),\bar a(\rho^\prime)\bigr)_q|0\rangle=
  \delta(\rho,\rho^\prime)~. \eqno(3.4)$$

\no Therefore in the absence of a background field, the usual
Feynman propagators are unchanged.  Under these conditions the
remaining effects of the $q$-commutators are relatively slight.
On the other hand it is easy to show by calculating spacelike field
commutators that this theory violates causality.$^3$

In a $q$-quantized theory, however, it is perhaps also natural to
consider the $q$-time ordered product

$$\eqalignno{T_q\bigl(\psi(x)\psi(x^\prime)\bigr) &=
  \psi(x)\psi(x^\prime) \quad t>t^\prime & (3.5a)\cr
  &=q\psi(x^\prime)\psi(x) \quad t<t^\prime & (3.5b)\cr}$$

\no This is obviously independent of the previous modification
of the canonical commutators.
When $q=-1$, (3.5) describes the usual $T$-product for Fermi
fields.

One may express the $q$-time ordered product as

$$T_q\bigl(\psi(x)\psi(x^\prime)\bigr)=
  {1\over 2}\biggl[\bigl\{\psi(x), \psi(x^\prime)\bigr\}_q+
  \epsilon(t-t^\prime)\bigl(\psi(x),\psi(x^\prime)\bigr)_q\biggr]
  \eqno(3.6)$$

$$\eqalign{\epsilon(t-t^\prime)&=1 \cr
  &=-1 \cr} \qquad
  \eqalign{& t>t^\prime \cr &t<t^\prime \cr} \eqno(3.7)$$

The naturalness of this product suggests that we examine the
$q$ modified $S$ matrix:

$$S^{(q)}=T_q\bigl(e^{i\int {\cal{L}}(x)d^4x}\bigr)~. \eqno(3.8)$$

\no In the remainder of this paper we shall discuss some features
of this quite different theory.
\vskip.5cm

\line {\bf 4. Normal Products. \hfil}
\s

To decompose (3.8) by Wick's theorem we begin by describing
normal products.

We first define the elementary normal products:

$$\eqalignno{N(1) &= 0 & (4.1)\cr
  N(\psi) &= \psi ~. & (4.2) \cr}$$

\no Then by (2.7)

$$N\bigl(a(\rho)\bar a(\rho^\prime)-
  q\bar a(\rho^\prime)a(\rho)\bigr)=
  N\delta(\rho,\rho^\prime)=0 \eqno(4.3)$$

\no and

$$\eqalignno{N\bigl(a(\rho)\bar a(\rho^\prime)\bigr)&=
  N\bigl(q\bar a(\rho^\prime)a(\rho)\bigr)=
  qN\bigl(\bar a(\rho^\prime)a(\rho)\bigr) & (4.4)\cr
  &=q\bar a(\rho^\prime)a(\rho)~. & (4.5)\cr}$$

The passage from (4.4) to (4.5) illustrates the general rule
that one must move all absorption operators to the right in order
to form a normal product.  At the same time one picks up a power
of $q$.

To get a string of absorption and emmission operators into normal
form there is no need to permute emission or absorption operators
among themselves; but every permutation of $a$ and $\bar a$ produces
$q$.  Therefore instead of the usual parity factor one now has
$q^n$.

One also has the usual decomposition into $+$ and $-$ frequency parts
since these are associated with emission and absorption operators.
Thus

$$\eqalignno{N(AB) &= N(A^++A^-, B^++B^-) \cr
  &= A^+B^++A^-B^-+A^+B^-+qB^+A^- \cr
  \noalign{\hbox{or}}
  N(AB) &= AB-(A^-,B^+)_q~ & (4.6)\cr
  \noalign{\hbox{and}}
  AB &=N(AB)+\langle AB\rangle & (4.7)\cr}$$

\no since $\langle AB\rangle$, the vacuum expectation value of $AB$,
is

$$\langle AB\rangle=\langle 0|A^-B^+|0\rangle
  = \langle 0|(A^-,B^+)_q|0\rangle  ~. \eqno(4.8)$$

Eq. (4.7) may be rewritten by (4.2)

$$N(A)B=N(AB)+\langle AB\rangle~. \eqno(4.9)$$

\no This relation may be generalized by induction to

$$N(A_1\ldots A_n)B=N(A_1\ldots A_nB) +
  \sum_{1\leq k\leq n} N(A_1\ldots
  \underbrace{A_k\ldots A_nB}) \eqno(4.10)$$

\no where $A_1\ldots A_nB$ are individually either absorption or
emission operators and

$$N(A_1\ldots \underbrace{A_k\ldots A_nB})=\eta\langle A_kB\rangle
  N(A_1\ldots A_{k-1}A_{k+1}\ldots A_n)~. \eqno(4.11)$$

\no Here $\eta=q^p$ results from the permutation of order between
the left and the right sides of (4.11).  $\langle A_kB\rangle$ is
called a pairing.

The proof of (4.10) follows along exactly the same lines as the
usual proof, since the only difference in the two situations
arises from the fact that usually $q=\pm 1$ and we finish
with $(\pm 1)^p$ while here we finish with $q^p$.

Wick's theorem for normal products follows immediately from (4.10)
and states that any product may be decomposed into a sum of normal
products with all possible pairings (including no pairings),
namely:

$$\eqalign{A_1\ldots A_n=N(A_1\ldots A_n) &+
  N(\underbrace{A_1A_2}\ldots A_n)
  +\ldots N(\underbrace{A_1\ldots A_{n-1}}A_n)
  +N(\underbrace{A_1\ldots A_n})\cr
  &+N(\underbrace{A_1A_2}\underbrace{A_3A_4}\ldots A_n)
  +\ldots~.\cr} \eqno(4.12)$$

\no This result is established by first assuming that it holds for
$n$, next multiplying by $A_{n+1}$ and using (4.10) to show that
it holds for $n+1$ and therefore by induction for all $n$.

To expand (3.8) one needs the $q$-chronological product (3.5) and
the corresponding pairing or vacuum expectation value

$$\overbrace{\psi_1(x_1)\psi_2}(x_2)=
  \langle 0|T_q(\psi_1(x_1)\psi_2(x_2)|0\rangle~. \eqno(4.13)$$

\no Now

$$T_q(\psi_1(x_1)\psi_2(x_2)=
  N(\psi_1(x_1)\psi_2(x_2)
  +\overbrace{\psi_1(x_1)\psi_2}(x_2) \eqno(4.14)$$

\no and

$$\eqalignno{\overbrace{\psi_1(x_1)\psi_2}(x_2) &=
  \underbrace{\psi_1(x_1)\psi_2}(x_2) \qquad
  x_1^o>x_2^o & (4.15a)\cr
  &= q\underbrace{\psi_2(x_2)\psi_1}(x_1) \qquad
  x_2^o>x_1^o & (4.15b)\cr}$$

\no Wick's theorem as applied to $q$-chronological products now
states that the $T_q$ product of a system of $n$ linear operators
is equal to the sum of their normal products with all possible
$q$-chronological pairings.  These pairings are the $q$-causal
propagators.
\vskip.5cm

\line {\bf 5. The $\Delta_q$ and $\Delta_{q_1}$ Functions. \hfil}
\s

By (3.6) the $q$-causal propagator is

$$\eqalignno{\Delta_{qF}(x-x^\prime) &=
  \langle 0|T_q\bigl(\psi(x),\psi(x^\prime)\bigr)|0\rangle & (5.1)\cr
  &= {1\over 2} \bigl[\Delta_1^q(x-x^\prime)+\epsilon(t-t^\prime)
  \Delta^q(x-x^\prime)\bigr] & (5.2)\cr}$$

\no where

$$\eqalignno{\Delta_1^q(x-x^\prime) &=
  \langle 0|\bigl\{\psi(x),\psi(x^\prime)\bigr\}_q|0\rangle & (5.3)\cr
  \noalign{\hbox{and}}
  \Delta^q(x-x^\prime) &=
  \langle 0|\bigl(\psi(x),\psi(x^\prime)\bigr)_q|0\rangle~. & (5.4)\cr}$$

\no In this section we shall calculate $\Delta^q$ and $\Delta^q_1$.

Let us adopt the representation (2.1).  Then to compute $\Delta^q$
and $\Delta^q_1$ we need

$$\eqalignno{\langle 0|\bigl(a(\rho),\bar a(\rho^\prime)\bigr)_q|0\rangle
  &= \langle 0|a(\rho)\bar a(\rho^\prime)-q~
  \bar a(\rho^\prime) a(\rho)|0\rangle\cr
  &= \delta(\rho,\rho^\prime) & (5.5)\cr
  \langle 0|\bigl(\bar b(\rho),b(\rho^\prime)\bigr)_q|0\rangle
  &= -q~\delta(\rho,\rho^\prime) & (5.6)\cr
  \langle 0|\bigl\{a(\rho),\bar a(\rho^\prime)\bigr\}_q|0\rangle
  &= \delta(\rho,\rho^\prime) & (5.7)\cr
  \langle 0|\bigl\{\bar b(\rho),b(\rho^\prime)\bigr\}_q|0\rangle
  &= q~\delta(\rho,\rho^\prime) & (5.8)\cr}$$

\no By (2.1)

$$\langle 0|\bigl(\psi_\alpha(x),\bar\psi_\beta(x^\prime)\bigr)_q
  |0\rangle=i\bigl(\Theta^{(+)}_{\alpha\beta}(p)
  \Delta_+(x-x^\prime)-q~\Theta^{(-)}_{\alpha\beta}
  \Delta_-(x-x^\prime)\bigr) \eqno(5.9)$$

\no where

$$\eqalignno{i~\Delta_\pm(x) &= \biggl({1\over 2\pi}\biggr)^3
  \int {d\vec p\over 2p_o}~e^{\mp i\rho x} & (5.10)\cr
  \Theta^{(+)}_{\alpha\beta}(\vec p) &= \sum_r
  u_\alpha(\vec p,r) \bar u_\beta(\vec p,r) & (5.11)\cr
  \Theta^{(-)}(\vec p) &= \sum_r v_\alpha(\vec p,r)
  \bar v_\beta(\vec p,r) & (5.12)\cr}$$

\no The $r$-sum is over spin states.  For scalars there is no spin
and therefore $\Theta^\pm=1$.  For spinors $u$ and $v$ refer to
positive and negative energy states.  (The bar means the Dirac
adjoint for spinors and the complex conjugate for scalars.)  Then
for scalars

$$\langle 0|\bigl(\psi(x),\bar\psi(x^\prime)\bigr)_q|0\rangle=
  i\Delta^-_q(x-x^\prime) \eqno(5.13)$$

\no where

$$\Delta^-_q(x)=\Delta_+(x)-q\Delta_-(x)~. \eqno(5.14)$$

\no For spinors

$$\Theta^\pm_{\alpha\beta}(\vec p)={1\over 2m}
  (\pm m+p \! \! \!/)_{\alpha\beta}~.  \eqno(5.15)$$

\no Then

$$\eqalignno{\langle 0|\bigl(\psi_\alpha(x),
  \bar\psi_\beta(x^\prime)\bigr)_q
  |0\rangle &= {1\over 2m} (m+i\partial \!\!\!/)_{\alpha\beta}
  \Delta_+(x)-q {1\over 2m}(-m-i\partial \!\!\!/)_{\alpha\beta}
  \Delta_-(x) \cr
  &= {1\over 2m}(m+i\partial \!\!\!/)_{\alpha\beta}
  \Delta^+_q(x) & (5.16)\cr}$$

\no where

$$\Delta^+_q(x)=\Delta_+(x)+q\Delta_-(x)~. \eqno(5.17)$$

\no But $q=\pm 1$ for scalars and spinors respectively.  Hence

$$\Delta_1^-(x)=\Delta^+_{-1}(x)~. \eqno(5.18)$$

\no Therefore in the limit, $|q|=1$, both fields are causal.

For neutral vector fields with no antiparticles we have

$$A_\alpha(x)=\biggl({1\over 2\pi}\biggr)^{3/2}
  \int {d\vec p\over (2p_o)^{1/2}} \sum_{r=1}^3
  \bigl[a(\vec p,r) e^{-ipx}+\bar a(\vec p,r)
  e^{ipx}\bigr] e_\alpha(\vec p,r) \eqno(5.19)$$

\no where $e_\alpha(\vec p,r)$ is the polarization vector.  Then

$$\langle 0|\bigl(A_\alpha(x),A_\beta(x^\prime)\bigr)_q|0\rangle
  =i~\Theta_{\alpha\beta}(\vec p)~\Delta^-_q(x-x^\prime) \eqno(5.20)$$

\no since

$$\langle 0|\bigl(\bar a(\rho),a(\rho^\prime)\bigr)_q
  |0\rangle =-q~\delta(\rho,\rho^\prime)~. \eqno(5.21)$$

\no Here

$$\Theta_{\alpha\beta}^{(\pm)}(\vec p)=
  \Theta_{\alpha\beta}(\vec p)=\sum_{r=1}^3
  e_\alpha(r,\vec p)e_\beta(r,\vec p) \eqno(5.22)$$

\no or

$$\Theta_{\alpha\beta}=g_{\alpha\beta}-{p_\alpha p_\beta\over
  m^2} \eqno(5.23)$$

\no and

$$\langle 0|\bigl(A_\alpha(x),A_\beta(x^\prime)\bigr)_q|0\rangle=
  i\biggl(g_{\alpha\beta}-{\partial_\alpha \partial_\beta\over m^2}\biggr)~
  \Delta^-_q(x-x^\prime) \eqno(5.24)$$

\no for a massive vector, while

$$\langle 0|\bigl(A^{\hbox{tr}}_i(x),A_j^{\hbox{tr}}(x^\prime)\bigr)_q
  |0\rangle=i\biggl(\delta_{ij}-{\partial_i\partial_j\over \partial^2}\biggr)~
  D^-_q(x-x^\prime) \eqno(5.25)$$

\no for the massless case.  Here

$$D(x)=\Delta(x;m=0) \eqno(5.26)$$

\no and

$$D^-_q=D^+-qD^-~. \eqno(5.27)$$

Similarly the field anticommutators are found to have the following
vacuum expectation values:
\s
\line {a) scalar \hfil}

$$\langle 0|\bigl\{\psi(x),\bar\psi(x^\prime)\bigr\}_q|0\rangle=
  \Delta^+_q(x-x^\prime) \eqno(5.28)$$

\line {b) vector \hfil}

$$\langle 0|\bigl\{A^{\hbox{tr}}_i(x),A^{\hbox{tr}}_j(x^\prime)\bigr\}_q
  |0\rangle=i\biggl(\delta_{ij}-
  {\partial_i\partial_j\over \partial^2}\biggr)~D^+_q(x-x^\prime) \eqno(5.29)$$

\no with a similar expression for the massive case.
\s
\line {c) spinor \hfil}

$$\langle 0|\bigl\{\psi_\alpha(x),\bar\psi_\beta(x^\prime)\bigr\}_q
  |0\rangle={1\over 2m}(m+i\partial \!\!\!/)_{\alpha\beta}~
  \Delta^-_q(x-x^\prime)~. \eqno(5.30)$$

In the boson examples the transition from commutator to anti-commutator
requires a change from $\Delta^-$ and $D^-$ to $\Delta^+$ and
$D^+$ while in the fermionic case the change is from $\Delta^+$ to
$\Delta^-$.
\vskip.5cm

\line {\bf 6. The Causal Propagator. \hfil}
\s

By (5.2) and the results of the previous section for $\Delta_q$
and $\Delta_{q1}$, the $q$-causal propagators in the three
cases are found to be the following:
\s
\line {a) scalar \hfil}

$$\Delta_{F_q}(x-x^\prime)={i\over 2}
  \bigl[\Delta_q^+(x-x^\prime)+\epsilon(t-t^\prime)
  \Delta_q^-(x-x^\prime)\bigr]  \eqno(6.1)$$

\line {b) massive vector \hfil}

$$\eqalign{(\Delta_{F_q})_{\alpha\beta}(x-x^\prime) &=
  \biggl(g_{\alpha\beta}-{\partial_\alpha\partial_\beta\over
  m^2}\biggr){i\over 2}\bigl[\Delta^+_q(x-x^\prime)+\epsilon(t-t^\prime)
  \Delta^-_q(x-x^\prime)\bigr]\cr
  &= \biggl(g_{\alpha\beta}-{\partial_\alpha\partial_\beta\over
  m^2}\biggr)~\Delta_{F_q}(x-x^\prime)\cr} \eqno(6.2)$$

\line {c) spinor \hfil}

$$\bigl(S_{F_q}(x-x^\prime)\bigr)_{\alpha\beta}={1\over 2m}
  (m+i\partial \!\!\!/)_{\alpha\beta}{i\over 2}
  \bigl[\Delta^-_q(x-x^\prime)+\epsilon(t-t^\prime)
  \Delta_q^+(x-x^\prime)\bigr]~. \eqno(6.3)$$
\s
\no In the boson examples we have

$$\eqalignno{\Delta_{F_q}&= {1\over 2}
  \bigl[(\Delta^{(+)}+q\Delta^{(-)}\bigr)+\epsilon(t)
  (\Delta^{(+)}-q\Delta^{(-)})\bigr]\cr
  &= {1\over 2}\bigl[(1+\epsilon(t))\Delta^{(+)}
  +q\bigl(1-\epsilon(t)\bigr)\Delta^{(-)}\bigr] & (6.4)\cr
  &= \Delta^{(+)} \quad t>0 \cr
  &= q\Delta^{(-)} \quad t<0  & (6.5)\cr}$$

\no or, by (6.10),

$$\eqalignno{i~\Delta_{F_q} &= \biggl({1\over 2\pi}\biggr)^3
  \int_H {d\vec k\over 2k_o}~e^{-ikx} \qquad t>0 & (6.6a)\cr
  &= q\biggl({1\over 2\pi}\biggr)^3 \int_H
  {d\vec k\over 2k_o}~e^{ikx} \qquad t<0 & (6.6b) \cr}$$

\no Here $H$ indicates integration over the mass hyperboloid

$$k_o=\omega=(\vec k^2+m^2)^{1/2}~. \eqno(6.7)$$

\no These expressions are equivalent to

$$\Delta_{F_q}(x) = \biggl({1\over 2\pi}\biggr)^4 \int_{F_q}
  {1\over k^2-m^2}~e^{-ikx} d^4k \eqno(6.8)$$

\no where the $F_q$ is the Feynman contour but the left hand pole
is of strength $q$, since (6.8) may be rewritten as follows:

$$\eqalign{\Delta_{F_q}(x)&=\biggl({1\over 2\pi}\biggr)^3 \int
  d\vec k~e^{i\vec k\vec x}{1\over 2\pi i} \int_{C_+}~
  e^{-ik_ot}{1\over k_o^2-\omega^2}~dk_o \quad t>0 \cr
  &=\biggl({1\over 2\pi}\biggr)^3 \int d\vec k~e^{i\vec k\vec x}
  {q\over 2\pi i}\int_{C_-}~e^{-ik_ot}{1\over k_o^2-\omega^2}~
  dk_o \quad t<0 \cr} \eqno(6.9)$$

\no Here $C_+$ and $C_-$ are clockwise contours in the complex
$k_o$ plane about the two points $\omega=\pm(k^2+m^2)^{1/2}$ on
the real axis.  The preceding equation (6.9) is equivalent to
(6.6).

An alternative way to write (6.8) is

$$\Delta_{F_q}(x)=\biggl({1\over 2\pi}\biggr)^4 \int_F
  \tilde\Delta_{F_q}(k)~e^{-ikx} d^4k \eqno(6.10)$$

\no where $F$ is again the Feynman contour and the Fourier
transform is

$$\eqalignno{\tilde\Delta_{F_q} &= {1\over 2\omega}
  \biggl({1\over k_o-\omega}-{q\over k_o+\omega}\biggr)  & (6.11)\cr
  &= {1\over 2}\biggl({1+q\over k^2-m^2}+
  {(1-q)(k_o/\omega)\over k^2-m^2}\biggr)~. & (6.12)\cr}$$

\no $\tilde\Delta_{F_q}$ is the propagator in momentum space for
the scalar field.

The corresponding propagators in momentum space for the vector
and spinor may be obtained from (6.2) and (6.3).  Except for spin
(6.1) and (6.2) are the same but (6.3) differs because there the
$\Delta^+_q$ and $\Delta^-_q$ functions are interchanged.
Therefore the spinor propagator in momentum space is obtained
from (6.12) and (6.3) by replacing $q$ by $-q$ as follows:

$${1\over 2m} (m+p\!\!\!/){1\over 2}
  \biggl({1-q\over p^2-m^2}+{(1+q)(p_o/\omega)\over
  p^2-m^2}\biggr)~. \eqno(6.13)$$

\no In the limit $|q|=1$ the scalar propagator becomes

$${1\over k^2-m^2} \eqno(6.14)$$

\no while the spinor propagator becomes

$${1\over 2m}(m+p\!\!\!/){1\over p^2-m^2}~. \eqno(6.15)$$

\no The $q$-modified propagators take the usual form not only
when $|q|=1$ but also when $|q|\not= 1$ and $k$ lies on the mass
hyperboloid $(k_o=\omega)$.  Therefore, if $|q|\not= 1$, the
dependence on $q$ becomes significant only for internal lines.
For the photon field we may write

$$D_{\mu\lambda}(x)=g_{\mu\lambda}\biggl({1\over 2\pi}\biggr)^4
  \biggl({1+q\over 2}\biggr) \int e^{-ikx}
  \biggl({1\over k^2-m^2}+{1-q\over 1+q}{k_o\over \omega}\biggr)
  d^4k \eqno(6.16)$$

\no and for the spinor field

$$S_{\alpha\beta}(x)=(i\partial\!\!\!/+m)_{\alpha\beta}
  \biggl({1\over 2\pi}\biggr)^4 \int e^{-ipx}
  {1\over 2}\biggl({1-q\over p^2-m^2}+
  {(1+q)(p_o/\omega)\over p^2-m^2}\biggr) d^4p~. \eqno(6.17)$$
\vskip.5cm

\line {\bf 7. Tests. \hfil}
\s

It is of course possible to examine the effects of these altered
propagators.  We consider two examples from QED:
\s
\line {a) Electron-Electron Scattering. \hfil}
\s

Let electrons in states $(A,B)$ scatter into states $(C,D)$:

$$A+B \rightarrow C+D~.$$

\no To lowest order the matrix element in configuration space is

$$\eqalign{\langle CD|S_q|AB\rangle &= {1\over 2}
  \biggl({e\over i\hbar c}\biggr)^2 \cr
   &\times \int \! \int
  d^4x_1d^4x_2\langle CD|N(J^\mu(x_1)J^\lambda(x_2))|AB\rangle
  \langle 0|T_qA_\mu(x_1)A_\lambda(x_2)|0\rangle \cr}\eqno(7.1)$$

\no where

$$\langle 0|T_q\bigl(A_\mu(x_1)A_\lambda(x_2)\bigr)|0\rangle=
  D^q_{\mu\lambda}(x_1,x_2) \eqno(7.2)$$

\no and

$$J^\mu(x)=\bar\psi(x)\gamma^\mu\psi(x) \eqno(7.3)$$

\no with

$$\psi(x)=\sum \bigl(a(\rho)u(\rho,x)+\bar b(\rho) v(\rho,x)\bigr)~.
  \eqno(7.4)$$

\no One finds

$$\langle CD|N(J^\mu(x_1)J^\lambda(x_2)\bigr)|AB\rangle=
  q(1+qP_{AB})(1+qP_{CD})
  \langle D|J^\mu(x_1)|B\rangle\langle C|J^\lambda(x_2)|A\rangle
  \eqno(7.5)$$

\no where

$$\langle D|J^\mu(x)|B\rangle=\bar u_D(x)\gamma^\mu u_B(x)~. \eqno(7.6)$$

\no Then

$$\eqalign{\langle CD|S_q|AB\rangle &=-i{q\over 2}
  \biggl({e\over \hbar c}\biggr)^2 \cr
  & \times\int \!\! \int d^4x_1d^4x_2(1+qP_{AB})(1+qP_{CD})
  \bigl(D|J^\mu(x_1)|B\bigr) \cr
  &\times \bigl(C|J^\lambda(x_2)|A\bigr)
  D^q_{\mu\lambda}(x_1,x_2)~.\cr} \eqno(7.7)$$

\no If $q=-1$ for the initial and final particles

$$D^{(1)}_{\mu\lambda}(x_1,x_2)=D^{(1)}_{\mu\lambda}
  (|x_1-x_2|)=D^{(1)}_{\mu\lambda}(x_2,x_1)~. \eqno(7.8)$$

\no Then

$$\eqalign{\langle CD|S_q|AB\rangle &=-i~q
  \biggl({e\over \hbar c}\biggr)^2 \cr
  & \times \int \!\! \int d^4x_1d^4x_2
  (1-P_{CD})\langle D|J^\mu(x_1)|B\rangle
  \langle C|J^\lambda(x_2)|A\rangle D^q_{\mu\lambda}(x_1-x_2)~.\cr}
  \eqno(7.9)$$

\no If $q\not= 1$,

$$D^q_{\mu\lambda}(x_1,x_2) \not= D^q_{\mu\lambda}(x_2,x_1)~. \eqno(7.10)$$

\no As a consequence there are four diagrams since the incoming as
well as the outgoing lines are crossed in order to describe the
$q$-antisymmetrization of the initial as well as the final state.

Let us consider just the contributions of the final state.  In
momentum space this is

$$\eqalign{\langle CD|S^2|AB\rangle \sim & q
  \bigl[\langle P_C|J^\mu(0)|P_A\rangle
  \Delta_{CA}\langle P_D|J_\mu(0)|P_B\rangle \cr
  &- \langle P_D|J^\mu(0)|P_A\rangle
  \Delta_{DA}\langle P_C|J_\mu(0)|P_B\rangle\bigr] \cr}\eqno(7.11)$$

\no where, by (6.16),

$$\eqalignno{\Delta_{CA}&={1\over |P_C-P_A|^2}~F_{CA} & (7.12)\cr
  \Delta_{DA} &= {1\over |P_B-P_A|^2}~F_{DA} & (7.13)\cr}$$

\no and

$$\eqalignno{F_{CA}&= {1+q\over 2}
  \biggl(1+{1-q\over 1+q}~{E_A-E_C\over |\vec P_A-\vec P_C|}\biggr)
  & (7.14)\cr
  F_{DA}&= {1+q\over 2}\biggl(1+{1-q\over 1+q}~
  {E_A-E_D\over |\vec P_A-\vec P_D|}\biggr) & (7.15)\cr}$$

\no If $q=1$, for the intermediate boson
the correction factors $F_{CA}$ and $F_{DA}$ are
not present and one gets the usual formula for Moller scattering.
If $q\not= 1$, these factors depend on frame.
In the center of mass system there is no effect.
\s
\line {b) Electron Positron Annihilation. \hfil}
\s

We consider

$$e_++e_- \rightarrow \gamma_1+\gamma_2~. \eqno(7.16)$$

\no The conservation of 4-momentum now reads

$$p_++p_-=k_1+k_2~. \eqno(7.17)$$

\no In the previous example there was an internal photon line.  In
this example there is an internal electron line.  By (6.17) the
usual matrix element for $e_+e_-$ annihilation is modified by the
same factors that appeared in the modified Moller formula, namely:

$$\eqalignno{F_{e_+k_1}&= \biggl({1-q\over 2}\biggr)
  \biggl(1+{1+q\over 1-q}{E_+-\epsilon_1\over
  |\vec p_+-\vec k_1|}\biggr) & (7.18)\cr
  F_{e_+k_2}&= \biggl({1-q\over 2}\biggr)
  \biggl(1+{1+q\over 1-q}
  {E_+-\epsilon_2\over |\vec p_+-\vec k_2|}\biggr)~. & (7.19)\cr}$$

\no In the center of mass system the numerators again vanish but
again the result is dependent on frame.  Therefore, although we have
assumed relativistic kinematics for the free particles, there remains
a non-relativistic frame dependence in the final result, coming from
the form of the internal propagators.  That is, the $q$-dynamics
explicitly violates relativity.
\vskip.5cm

\line {\bf 8.  Remarks. \hfil}
\s

The preceding work has been based on two changes in the postulational
basis of QED: first the use of $q$-commutators at equal times and
second, the use of $q$-causal propagators.  At weak field strengths
or in the absence of background fields the use of $q$-commutators
at equal times ordinarily makes very little difference if $q$ is sufficiently
close to unity; but the theory is still non-causal.
On the other hand, the use of $q$-causal propagators
leads to non-relativistic results and thus merges with an alternative
approach in which the Poincar\'e group is explicitly deformed.  In
both cases deformation of dynamical law is not independent of
deformation of the geometry.
\vskip.5cm

\line {\bf References. \hfil}
\s
\item{1.} E. Marcus and R. Finkelstein, J. Math. Phys. (in press);
UCLA/94/TEP/32.
\item{2.} O. Ogievetsky, W. B. Shmidke, J. Wess, and B. Zumino,
Lett. Math. Phys. {\bf 23}, 233 (1991).
\item{3.} R. Finkelstein, Lett. Math. Phys. (in press); UCLA/94/TEP/26.
\ve

\bye